\title{The problem of prescribed critical functions}
\author{Emmanuel Humbert\thanks{ Institut \'Elie Cartan, Universit\'e de Nancy 
  1, BP 239 54506 Vandoeuvre-L\`es-Nancy Cedex, FRANCE  
Email: humbert@iecn.u-nancy.fr  }  
 and Michel Vaugon\thanks{   
  Institut de math\'ematiques de Jussieu, Universit\'e de Paris 6, 175, rue 
  du Chevaleret 75013 Paris, France - Email: vaugon@math.jussieu.fr }} 
\documentclass[12pt]{article} 
\usepackage{amsmath} 
\numberwithin{equation}{section} 
\newtheorem{Defi}{Definition}[section] 
\newtheorem{theorem}{Theorem}[section]

\newtheorem{step}{Step} 
 
\begin{document} 
 
%%%%%%%%%%%%%%%%%%%%%%%%%%%%%%%%%%%%%%%%%%%%%%%%%%%%%%%%%%%%%%%% 
\maketitle 
%%%%%%%%%%%%%%%%%%%%%%%%%%%%%%%%%%%%%%%%%%%%%%%%%%%%%%%%%%%%%%%%nash 
 
\begin{abstract} 
Let $(M,g)$ be a compact Riemannian manifold on dimension $n \geq 4$ not 
conformally diffeomorphic to the sphere $S^n$. We 
prove that a smooth function $f$ on $M$ is a critical  
function for a metric $\tilde{g}$ conformal to 
$g$ if and only if there exists $x \in M$ such that $f(x)>0$. 
%
%\begin{center}
%{\bf R\'esum\'e}
%\end{center}
% 
%Soit $(M,g)$ une vari\'et\'e Riemannienne compacte de dimension $n \geq 4$
%non conform\'ement diff\'eomorphe \`a la sp\`ere $S^n$. Nous montrons
%qu'une fonction lisse $f$ sur $M$ est critique pour une m\'etrique
%$\%tilde{g}$ conforme \`a $g$ si et seulement si il existe un point $x \in
%M$ pour lequel $f(x) > 0$.
\end{abstract} 
%%%%%%%%%%%%%%%%%%%%%%%%%%%%%%%%%%%%%%%%%%%%%%%%%%%%%%%%%%%%%%%% 
{\bf Keywords:}          
Best constants, Sobolev inequalities.          

\noindent {\bf Mathematics Classification:}       
53C21, 46E35, 26D10.

\section{Introduction} 
\subsection{Critical functions}

Let $(M,g)$ be a compact Riemannian manifold on dimension $n \geq 3$. The 
Sobolev embedding $H_1^2(M)$ into $L^N(M)$ ($N=\frac{2n}{n-2}$) asserts 
that there exists two constants $A,B>0$ such that, for all $u \in H_1^2(M)$,  
$${\left( \int_M {|u|}^N dv_g \right)}^{\frac2N} \leq  A \int_M {|\nabla u|}^2 dv_g +  
\int_M B u^2 dv_g \eqno{S(A,B)}$$ 
Here, $H_1^2(M)$ is the set of functions $u \in L^2(M)$ such that $\nabla u 
\in L^2(M)$. It is well known that the best constant $A$ in 
  this inequality is 
$$A= K(n,2)^2= \frac{4}{n(n-2)\omega_n^{\frac{2}{n}}}$$ 
where $\omega_n$ stands  for the volume of the standard $n$-dimensional
sphere. As shown by Hebey and Vaugon \cite{hv1}, this best constant is
attained. In  
other words, there exists $B>0$ such that $S(K(n,2)^2,B)$ is true for all $u 
\in H_1^2(M)$. We note $B_0(g)$ the smallest constant $B$ such this 
assertion is true. Clearly,  $S(K(n,2)^2,B_0(g))$ 
holds for all $u \in H_1^2(M)$. This inequality is sharp. 
Moreover, we have (see the general reference \cite{heb})
$$B_0(g) \geq \max \big(  
\frac{n-2}{4(n-1)} K(n,2)^2 \max_{M} S_g, Vol_g(M)^{-\frac{2}{n}} \big)  
\eqno{(*)}$$ 
where $S_g$ the scalar curvature of $g$.  
Natural questions are then: \\

\noindent -does there exists extremal functions in 
$S(K(n,2)^2,B_0(g))$? 
Extremal functions are nonzero functions for which  
$S(K(n,2)^2,B_0(g))$ is an equality.  
 
\noindent -is it possible that $(*)$ is an equality?

\noindent These questions seem to be independent but Djadli and Druet 
 proved in \cite{dd} that one of the following assertions must hold if $n
 \geq 4$: \\

\noindent a)$B_0(g) =  \frac{n-2}{4(n-1)}  
K(n,2)^2 \max_{M} S_g$; 
 
\noindent b)$S(K(n,2)^2,B_0(g))$ possesses extremal 
functions. \\

\noindent Then, other questions arises naturally:  
 
\noindent -is it possible that a) is true and b) is false? 

\noindent -is it possible that b) is true and a) is false? 

\noindent -is it possible that a) and b) are true?

\noindent  Critical functions have been  introduced by Hebey and Vaugon in 
\cite{hv} in the purpose of answering this type of questions. The idea was to  
consider inequality $S(K(n,2)^2,B_0(g))$ in a metric $\tilde{g}$ conformal
to $g$. Namely, if $\tilde{g}= u^{\frac{4}{n-2}} g$ where $u \in
C^{\infty}(M)$, $u>0$, then one may check that
inequality $S(K(n,2)^2,B_0(g))$ is equivalent to the following one: for all
$u \in H_1^2(M)$, we have 
$${\left( \int_M {|u|}^N dv_{\tilde{g}} \right)}^{\frac2N} \leq  K(n,2)^2 
\int_M {|\nabla u|}^2 dv_{\tilde{g}} +  
\int_M f u^2 dv_{\tilde{g}} \eqno{S'(f,\tilde{g})}$$ 
where $f \in C^{\infty}(M)$ satisfies
$$\Delta_g u + B_0(g) u = f u^{N-1}$$
Note that this implies that $B_0(\tilde{g}) \leq \max(f)$.
It is then natural to introduce the notion of critical function. 
Critical functions corresponds to ``best functions'' in inequality
above. More precisely, 

\begin{Defi} [Hebey, Vaugon \cite{hv}] 
We say that a smooth function $f$ is critical for a metric $g$ if 
$S'(K(n,2)^2 f,g)$ is true for all $u \in H_1^2(M)$ and if for all smooth
function $f' \leq f$ with $f' \not= f$, inequality $S'(K(n,2)^2 f',g)$ is not
true.
\end{Defi}
 
\noindent Another way to define critical functions is the following.
For any $u \in 
H^2_1(M)-\{0\}$, we define: 
$$I_{\tilde{g},f}(u)= \frac{\int_M {|\nabla u|}^2_{\tilde{g}} dv_{\tilde{g}} +  
\int_M f u^2 dv_{\tilde{g}}}{  {\left( \int_M {|u|}^N 
dv_{\tilde{g}} \right) }^{\frac2N} }   $$ 
and  
$$\mu_{\tilde{g},f}= \inf_{u \in H_1^2(M)-\{0\}} I_{\tilde{g},f}(u)$$ 
It is well known that  
$$\mu_{\tilde{g},f} \leq  K(n,2)^{-2}$$ 
We now say that
 
\begin{Defi}  [Hebey, Vaugon \cite{hv}] 
A smooth function $f$ is   
\begin{itemize}  
\item subcritical for $\tilde{g}$ if $\mu_{\tilde{g},f} < K(n,2)^{-2}$; 
\item weakly critical for $\tilde{g}$ if  
$\mu_{\tilde{g},f} = K(n,2)^{-2}$; 
\item critical for  $\tilde{g}$ if $f$ is weakly critical with the 
  additional property that  
  for any smooth function $h$  such that $h \leq f$ 
and $h \not= f$,  $h$ is subcritical. 
\end{itemize} 
 
\end{Defi}

From these definitions, we get some remarks. At first,
let $\tilde{g}$ be a metric conformal to $g$ and $f,h$ 
two smooth functions on $M$ such that $h \leq f$. Then, it is clear that:
for all $u \in 
H^2_1(M)-\{0\} $, $I_{\tilde{g},f}(u) \leq I_{\tilde{g},h}(u)$. Hence, the
fact that  
$h$ is weakly critical for $\tilde{g}$ implies that  $f$ is  weakly
critical for $\tilde{g}$ and the fact that $f$ is subcritical for $\tilde{g}$
implies $h$ is  subcritical for $\tilde{g}$.\\ 
 
\noindent A second remark is that a weakly  
critical function $f$ for $\tilde{ g}$ satisfies: 
$$f \geq \frac{n-2}{4(n-1)}S_{\tilde{g}}$$
As one can check, this can be proved by  mimicking the proof of $(*)$ in
\cite{heb}.\\ 
 
\noindent 
Assume now that $f$ is weakly 
critical for $\tilde{g}$ and that there exists $u \in 
H_1^2(M)-\{0\}$ such  
that $I_{\tilde{g},f}(u)= K(n,2)^{-2}$. Then, $f$ is critical for 
$\tilde{g}$. Indeed, if $h \leq f$ and $h \not= f$, we have  
 $I_{\tilde{g},h}(u) < I_{\tilde{g},f}(u)= K(n,2)^{-2}$ and hence, $h$ is 
 subcritical. A first consequence of this remark is that if $n \geq 4$, if $g$ is a metric
 such that $S_g$ is a constant function and if $M$ is not conformally
 diffeomorphic to the standard sphere, then $ K(n,2)^{-2} B_0(g)$ is a
 critical function for $g$. Indeed, since $M$ is not conformally
 diffeomorphic to the standard sphere, it is well known 
$\frac{n-2}{4(n-1)} S_g$ is subcritical for $g$. Hence,  
$B_0(g)> \frac{n-2}{4(n-1)} S_gK(n,2)^2$. By Djadli and
Druet's work \cite{dd},  there exists $u \in 
H_1^2(M)-\{0\}$ such that $I_{g,  K(n,2)^{-2} B_0(g)}(u)= K(n,2)^{-2}$. The
remark above gives then the result.

\noindent A second consequence is that this gives a third
definition of critical functions. Namely, let $f$ be a critical function
for $g$. Then, for any $f' \leq f$, $f' \not= f$, we have $\mu_{g,f'} <
K(n,2)^{-2}$. It is well known that this implies that $\mu_{g,f'}$ is
attained. In other words, there exists a minimizing solution of equation 
$$\Delta_g u + f' u= u^{N-1} \eqno{(**)}$$
Moreover, by remark above, for any smooth function $f' \geq f$, $f' \not=
f$, $\mu_{g,f'}$ is not attained and hence, equation (**) does not possess
any minimizing solution. Reciprocally, let $f$ be a smooth function which
satisfies these properties. Then, clearly $f$ is weakly critical. Hence,
there exists a critical fonction $h \leq f$ (see \cite{hv}). Assume that $h
\not=f$. Then, if $h'$ is a smooth function such that $h \leq h' \leq f$,
we know that equation (**) cannot have minimizing solution. This
contradicts the definition of $f$ and hence, $h=f$ and $f$ is a critical
for $g$. 
We have proven that the following definition is equivalent to the two
definitions given above.

\begin{Defi}
A smooth function $f$ is critical for a metric $g$ if it satisfies the two
following properties:

\noindent - for any  $f' \leq f$, $f' \not=
f$, equation (**) has a minimizing solution;

\noindent - for any  $f' \geq f$, $f' \not=
f$, equation (**) does not have minimizing solutions.
\end{Defi}

A very important property of critical functions is that they have a ``good''
transformation law when we make a conformal change of metric. Indeed, we set
$\tilde{g}= u^{\frac{4}{n-2}} g$ where
$u\in C^{\infty}(M)$, $u>0$.  
Let  $h \in C^{\infty}(M)$ and $v \in H^2_1(M)-\{0\}$. 
Then, standard computations show that 
$$I_{g,h}(v)= I_{\tilde{g}, \tilde{h}}(u^{-1} v)$$ 
where $\tilde{h}$ and $h$ are related by the following equation: 
 
$$\Delta_g u + h u = \tilde{h} u^{\frac{n+2}{n-2}} $$ 
This implies that $h$ is critical for $g$ if and only if $\tilde{h}= 
\frac{\Delta_g u + h u}{u^{\frac{n+2}{n-2}}}$ is critical for $\tilde{g}$. 
Another way to present this result is to say: 
if $f$ is a smooth function, then $f$ is critical for   
$\tilde{g}= u^{\frac{4}{n-2}} g$ if and only if $f u^{\frac{4}{n-2}}-  
\frac{\Delta_g{u}}{u}$ is critical for $g$.\\ 

%%%%%%%%%%%%%%%%%%%%%%%%%%%%%%%%%%%%%%%%%%%%%%%%%%%%%%%%%%%%%%%%%%%%% 

\subsection{The problem}

\noindent  In this paper, we consider the problem of prescribed critical function: 
let $(M,g)$ be a compact Riemannian manifold of dimension $n \geq 4$ not 
conformally diffeomorphic to the sphere $S^n$ 
and $f$ be a smooth function. Does there exist a metric $\tilde{g}$  
conformal to $g$ such that $f$ is a 
critical function for $\tilde{g}$? 
As explained in section 1,
critical functions plays an important role in the study of sharp Sobolev
inequalities. Therefore, critical functions must be studied deeply to
understand better these inequalities. Moreover, this problem is closely
related to important geometric problems as Yamabe problem or prescribed scalar
curvature problem. Namely, for two functions $\alpha, \beta \in C^{\infty}$, we
consider the following equation:
$$\Delta_g u+ \alpha  u = \beta u^{N-1}$$
This type of equation is very important in geometry. For example, the
Yamabe problem (see \cite{a}) consists in finding a smooth strictly positive
solution $u$ 
of this equation where $\alpha= \frac{n-2}{4(n-1)} S_g$ and where $\beta$ is
a constant function. On what concerns the prescribed scalar curvature
problem (see again \cite{a}), 
we have to find a smooth strictly positive solution $u$ where 
$\alpha= \frac{n-2}{4(n-1)} S_g$ and where $\beta$ is given. In our problem, it
follows from last paragraph that we are lead to find a critical function
$h$ for the metric $g$ and   
a smooth strictly positive
solution $u$  
of this equation where $\alpha=h$  and where $\beta=f$. 
Then, setting $\tilde{g}=u^{\frac{4}{n-2}}g$, we
obtain a conformal metric to $g$ for which $f$ is critical.
As in the prescribed scalar curvature problem, the difficulty here 
comes from the fact that it  
cannot be solved by variationnal methods. We give a complete resolution of
the problem in dimension $n \geq 4$. This is object of the following result:

\noindent {\bf Main theorem} 
{\it Let $(M,g)$ be a compact Riemannian manifold of 
dimension  $n \geq 4$ not conformally diffeomorphic to the sphere $S^n$ and 
let $f$ a smooth function on $M$. Then, 
there exists a metric $\tilde{g}$ conformal to $g$ for which $f$ is critical 
if and only  if there  
exists $x \in M$ such that $f(x)>0$. }\\ 
 
\noindent Obviously, the difficulty is to show that if  there  
exists $x \in M$ such that $f(x)>0$ then we can find a metric $\tilde{g}$ 
conformal to  $g$ for which $f$ is critical. Moreover,  
one can find a proof much more 
easier that the one we give here when $f > 0$. 
In other words, the difficult  
part of the theorem 
corresponds to the case when $f$ changes of sign. 
   
\noindent One can consider the same problem if $(M,g)$ is 
conformally diffeomorphic to the standard sphere or if $n=3$. At first, let 
$(S_n,h_0)$ be the standard sphere of dimension 
$n$ and $g$ be a metric conformal to $h_0$. Then the only critical function 
for $g$ is $\frac{n-2}{4(n-1)}S_g$. Hence, the problem is equivalent to  
the problem of prescribed scalar curvature. If now $n=3$, then critical 
functions do not have the same  properties than in dimension upper than four.  
For example, theorem 2.1 below is false when $n=3$ (see \cite{hv} and
\cite{d}). 
The 
$3$-dimensional  
case seems to be  interesting but the methods used here are not adapted to 
this case. \\

\section{Proof of main theorem} 
In this section, $(M,g)$ is  a compact Riemannian manifold of 
dimension  $n \geq 4$ not conformally diffeomorphic to the sphere $S^n$ and 
$f$ is a smooth function on $M$. In addition, up to making a conformal 
change of metric, one may assume that $S_g$ is a constant function and up to
multiplying $g$ by a constant, we can also assume that

\begin{eqnarray} \label{d1}
\alpha_0  - \max_M(f) \geq {n-2 \over 4(n-1)} S_g
\end{eqnarray}

\noindent where $\alpha_0 = K(n,2)^{-2} B_0(g)$.

\subsection{A preliminary result} 
For the proof of main theorem, we will need the following result: 
 
\begin{theorem} 
Let $(h_m)_m$ be a sequence of smooth functions on $M$ which converges 
uniformly to a smooth function $h$. We assume that for all $m$, $h_m$ is 
subcritical for $g$ and that $h$ is weakly critical. Moreover, we assume 
that  
$$h > \frac{n-2}{4(n-1)}S_g$$ 
Then, $h$ is critical. 
\end{theorem} 
The proof follows very closely the  
proof of Druet and Djadli's theorem  in 
\cite{dd}. In addition, the reader may refer to \cite{hv} sections 2 and 3 
for a sketch of  
 proof of this theorem as stated here.\\ 
  
% \noindent {\bf Sketch of proof:} 
 
%\noindent At first, since $h_m$ is subcritical, standard methods (see for 
%example \cite{a} or \cite{h}) show that 
%there exists $u_m \in \mathcal{H}_g$, $u_m >0$ such that  
%$$I_{g,h_m}(u_m)=\mu_{g_{h_m}} < K(n,2)^{-2}$$ 
%Writing Euler equation for $u_m$, we get that 
%$$\Delta_g u_m + h_m u_m = \mu_m u_m^{\frac{(n+2)}{(n-2)}}  
%\eqno{(E_m)}$$ 
%where, for convenience, we have written  
%$$\mu_m= \mu_{g_{h_m}}$$ 
%Clearly, the sequence $(u_m)_m$ is bounded in $H_1^2(M)$. By standard 
%methods, there exists $u \in H_1^2(M)$, $u \geq 0$ such that such that, up 
%to a subsequence, $(u_m)_m$ tends to $u$ weakly in $H_1^2(M)$, strongly in 
%$L^2(M)$ and almost everywhere. Passing to the limit in $m$, we get from 
%equation $(E_m)$ that $u$ is a weak solution of 
%$$\Delta_g u + h u = \mu u^{\frac{(n+2)}{(n-2)}} $$ 
%where $\mu=\lim_m \mu_m \leq K(n,2)^{-2}$ which exists up to a 
 %subsequence. Assume   
%that $u \not\equiv 0$. Then, by regularity theorems and maximum 
%principle, $u$ is smooth and positive. As one can check, 
%we have $u \in \mathcal{H}_g$. Clearly, this implies that 
%$$I_{g,h}(u)= \mu \leq  K(n,2)^{-2}$$ 
%Since $h$ is weakly critical then we get 
%$$I_{g,h}(u)=   K(n,2)^{-2}$$ 
%By remark 2 of section 1, $h$ is critical. 
%Now, we assume that $u \equiv 0$. 
%By standard arguments, there exists $x_0 \in M$ such that  for all $\delta 
%>0$,  
%$$\lim_m \int_{B(x_0,\delta)} u_m^N dv_g =1$$ 
%Moreover, we have: 

\subsection{Proof of main theorem} 
At first, if $f$ is a critical function for any metric $\tilde{g}$ (not
necessarily conformal to $g$), then, there exists $y \in M$ such that
$f(y)>0$. Indeed, coming back to the notations of section 1, if $f \leq 0$,
we have $I_{\tilde{g},f}(1) \leq 0 < K(n,2)^{-2}$ and hence, $f$ cannot be
critical for $\tilde{g}$. Therefore, we assume that   there exists $y \in M$ 
such that
$f(y)>0$ and we have to show that we can find a metric $\tilde{g}$
conformal to $g$ for which $f$ is critical. We set 
\[ \mathcal{F} :\left| \begin{array}{ccc} 
\Omega  & \to & C^{\infty}(M) \\ 
u & \mapsto &fu^{\frac{4}{n-2}} - \frac{\Delta_g u}{u} 
\end{array} \right.\] 
where  
$$\Omega=\{ u \in C^{\infty}(M) | u>0 \}$$ 
Let $u \in \Omega$. By paragraph 1.1, $f$ is critical for 
$\tilde{g}=u^{\frac{4}{n-2}}g$ if and only if $\mathcal{F}(u)$ is critical 
for $g$.  
 In the following, we say weakly critical, subcritical and critical
and we omit to say ``for $g$''. Coming back to the notations of section 1,
we set,  
for any smooth function $h$ and $u \in H_1^2(M)$: 
$$I_h(u) = I_{g,h}(u)= \frac{\int_M {|\nabla u|}^2_g dv_g+  \int_M h  u^2 dv_g} 
{{\left( \int_M {|u|}^N dv_g \right)}^{\frac{2}{N}}} $$ 
and  
$$\mu_h= \mu_{g,h}=\inf_{H_1^2(M)-\{0\}} I_h$$ 
With these notations, we want to find $u  \in \Omega$ such that 
$\mu_{\mathcal{F}(u)} = K(n,2)^{-2}$ and for all $h \leq \mathcal{F}(u)$, 
$h \not= \mathcal{F}(u)$, $\mu_h <K(n,2)^{-2}$. 
  A natural idea to prove the main theorem is then the following:   
we find $u,v 
\in \Omega$ such that $\mathcal{F}(u)$ is subcritical  and 
$\mathcal{F}(v)$ is weakly critical. Then, we take a continuous path  
$(u_t)_t \subset \Omega$ for $t \in [0,1]$ such that $u_0=u$ and $u_1=v$.  
We define 
$$t_0= \inf \{ t>0 \hbox{ such that } \mathcal{F}(u_t) \hbox{ is weakly critical } \}$$  
The idea is  to apply theorem 2.1 with $h=\mathcal{F}(u_{t_0})$ and  
$h_m=\mathcal{F}(u_{t_m})$ where 
$t_m=t_0-\frac1m$. The difficulty is that we 
need the additional assumption that $\mathcal{F}(u_{t_0}) > 
\frac{n-2}{4(n-1)}S_g$. The linear transformation $u_t = 
tu+(1-t)v$ does not work in general. Hence, we must be very careful with 
the choice of $u$, $v$ and $u_t$. In fact, we show that there exists $u \in
\Omega$ and $s>1$  
such that $\mathcal{F}(u)$ is subcritical and such that  $\mathcal{F}(u^s)$ is 
weakly critical. The method described above then works with $u_t= 
\mathcal{F}(u^{ts})$. We strongly use in the whole proof concentration phenomenoms. 
In the special case where $f > 0$ then,  
one can find a shorter proof than 
the one we give here. \\ 
 
\noindent Let us start the proof now. For all $u \in H^2_1(M)$, $t>0$, and  
 $q \in ]2, N]$, we set  
$$J_{q,t}(u)= \frac{\int_M {|\nabla u|}^2_g dv_g+  t \int_M   u^2 dv_g} 
{{\left( \int_M f {|u|}^q dv_g \right)}^{\frac{2}{q}}}$$ 
We define, for $q \in ]2, N]$ 
$$\mu_{q,t}= \inf_{\mathcal{H}_q} J_{q,t}(u)$$ 
where  
$$ \mathcal{H}_q =\{ u \in H^2_1(M) |  
\int_M f |u|^q dv_g  > 0 \}$$ 
Obviously, since there exists $y \in M$ such that $f(y)>0$ and since $f$ is 
continuous, the set  
$ \mathcal{H}_q$ is not empty. 
It is well known that for all $t >0$ 
\begin{eqnarray} \label{r0}
\mu_{N,t} \leq K(n,2)^{-2} {(\max_M{f})}^{\frac2N} 
\end{eqnarray}
and that for all $u \in  \mathcal{H}_q $, we have 
$J_t (u)=J_t(|u|)$. 
Hence, we can replace $  \mathcal{H}_q $ by  
 
$$ \mathcal{H}_q =\{ u \in C^{\infty}(M) | u>0 \hbox{ and } 
\int_M f u^q dv_g >0 \}$$

\noindent We now define: 
 
$$\Omega_{q,t}= \{ u \in \mathcal{H}_q | J_{q,t}(u)=\mu_{q,t} \hbox{ and }  
 \int_M f {u}^q dv_g = \mu_{q,t}^{\frac{q}{q-2}} \}$$ 
 
\noindent The value $ \mu_{q,t}^{\frac{q}{q-2}}$ is chosen to obtain 
equation $E(q,t)$ below. Note that, for all $t>0$ and $q \in ]2,N[$,  
$ \mu_{q,t} > 0$.  
By standard elliptic theory, we know that for all $q < N$ and all 
$t>0$, 
$$\Omega_{q,t} \not= \emptyset$$ 
Note that if $t$ is large (for example $t> \alpha_0=B_0(g) K(n,2)^{-2}$), 
$\Omega_{N,t} = \emptyset$. 
Indeed, let $t>\alpha_0$ and assume that there exists  
$u \in \Omega_{N,t}$.  
Then, by (\ref{r0}),  
$$ {(\max_M{f})}^{-\frac2N} \frac{\int_M {|\nabla u |}^2 dv_g + \alpha_0 \int_M u^2}{{\left( \int_M u^N dv_g 
    \right)}^{\frac{2}{N}} } <  J_{N,t}(u)$$ 
$$= \mu_{N,t} \leq K(n,2)^{-2} 
{(\max_M{f})}^{-\frac2N}$$  
This contradicts the fact that, by definition $\alpha_0$ is weakly
critical. 
Another remark is the following: if $u \in \Omega_{q,t}$ with $q \in ]2,N]$ 
and $t>0$, then writing the 
Euler equation of $u$, we get that $u$ satisfies  
$$\Delta_g u + t u= f u^{q-1}  \eqno{E(q,t)}$$

\noindent We first prove that: 
 
\begin{step} \label{s1} 
Let $0< t< \alpha_0$ where $\alpha_0= B_0(g) K(n,2)^{-2}$ is the lowest weakly
critical 
constant function.  
Then, there exists $q_0 <N$ such that for all $q \in 
[q_0,N[$, and all $u \in \Omega_{q,t}$, $\mathcal{F}(u)$ is subcritical  
for $g$. 
\end{step} 
  
\noindent { 
We proceed by contradiction. We assume that there exits a sequence $(q_i)$ 
of real numbers and a sequence $(u_i)$ of functions belonging to  
$\Omega_{q_i,t}$  such that   
 
\noindent - $\lim_i q_i=N$ 
 
\noindent - $q_i < N$ for all $i$ 
 
\noindent - $\mathcal{F}(u_i)$ is weakly critical.\\ 
 
\noindent  Clearly, $(u_i)$ is bounded in $H_1^2(M)$. Hence, by
standard arguments 
(see \cite{a} or \cite{h}),  
there exists  $u \in H_1^2(M)$ such that, up to a subsequence,  
$ u_i \to u$ weakly in   $ H_1^2(M)$, strongly in $L^2(M)$ 
, strongly in  $L^{N-2}(M)$  
and almost everywhere.  
  
\noindent First, we assume that $u \not\equiv 0$. Then, by elliptic theory, 
$u \in \mathcal{H}_N$ and up to subsequence, we may assume that  
$$ u_i \to u \hbox{ in }  C^2(M)$$ 
Therefore, the sequence $(\mathcal{F}(u_i))$  converges uniformly to 
$\mathcal{F}(u)$. Since $\mathcal{F}(u_i)$ is weakly critical, then 
$\mathcal{F}(u)$ is weakly critical too. Moreover, $u_i \in \Omega_{q_i,t}$
and hence  
satisfies equation $E(q_i,t)$. This gives 
 
\begin{eqnarray} \label{r2} 
\mathcal{F}(u_i) = t+ f(u_i^{N-2} -u_i^{q_i-2}) 
\end{eqnarray} 
 
\noindent Passing to the limit in $i$,  
we get that  
$$\Delta_g u + t u = f u^{N-1}$$ 
and 
$$\mathcal{F}(u) = t$$ 
Therefore, we have proven that the constant function  
$t$ is weakly critical for $g$.  
This is impossible 
since $t< \alpha_0$ and since $\alpha_0$ is the smallest weakly critical
constant function.\\

\noindent We now deal with the case where $u \equiv 0$. Since $t < 
\alpha_0$, the constant function $t$ is subcritical. Hence, there exists a 
positive function  
$\phi \in C^{\infty}(M)$ such that  
 
\begin{eqnarray} \label{r1} 
I_t(\phi)= \frac{\int_M {|\nabla \phi|}^2_g dv_g+ t \int_M  \phi^2 dv_g} 
{{\left( \int_M {|\phi|}^N dv_g \right)}^{\frac{2}{N}}}   
 < K(n,2)^{-2}  
\end{eqnarray}

\noindent Plugging the test function $\phi$ into $I_{\mathcal{F}(u_i)}$, we get 
by (\ref{r2}) that  
$$I_{\mathcal{F}(u_i)}(\phi)= I_t(\phi) + \frac{\int_M f( u_i^{N-2} -u_i^{q_i-2}) \phi^2 dv_g}{{\left( \int_M {|\phi|}^N dv_g \right)}^{\frac{2}{N}}}$$ 
By strong convergence of $(u_i)$ to $0$  
in $L^{N-2}(M)$ and since $q_i - 2 \leq N-2$, 
 we get that 
$$\lim_i \int_M f u_i^{N-2} \phi^2 dv_g  = \lim_i \int_M f u_i^{q_i-2}
\phi^2 dv_g  
=0$$ 
It follows that  
$$\lim_i I_{\mathcal{F}(u_i)}(\phi) =  I_t(\phi)  < K(n,2)^{-2} $$ 
That contradicts the fact that $\mathcal{F}(u_i)$ is weakly critical. 
This ends the proof of step \ref{s1}. \\ 
 
\noindent By step \ref{s1}, one can find two sequences $(q_i)$, $(t_i)$  
of real numbers such that  
 
\noindent 1) $2 <q_i < N$ and $\lim_i q_i= N$; 
 
\noindent 2) $t_i > 0$ and $\lim_i t_i =\alpha_0= K(n,2)^{-2} B_0(g)$ 
 
\noindent and a sequence $(v_i)$ of functions belonging to $\Omega_{q_i,t_i}$ 
 with the  
 additionnal property that $\mathcal{F}(v_i)$ is subcritical. 
Clearly, proceeding as in step \ref{s1}, one can find  $v \in H_1^2(M)$ such that,  
up to a subsequence,  
$ v_i \to v$ weakly in   $ H_1^2(M)$, strongly in $L^2(M)$ 
, strongly in  $L^{N-2}(M)$  
and almost everywhere. We set 
$$J_i=J_{q_i,t_i} \hbox{ and } \mu_i=\mu_{q_i,t_i}$$
We prove that

\begin{step} \label{s2} 
We can assume  $v \equiv 0$ 
\end{step}  
Otherwise, as in step \ref{s1},  
$$v_i \to v \hbox{ in } C^2(M)$$ 
Mooreover, $v_i \in \Omega_{q_i,t_i}$. Hence, $v_i$ satisfies equation 
$E(q_i,t_i)$ and we have  
 
\begin{eqnarray} \label{r3} 
\mathcal{F}(v_i) =  t_i + f ( v_i^{N-2}- v_i^{q_i-2})  
\end{eqnarray} 
 
\noindent Passing to the limit in $i$, we get  
$$\mathcal{F}(v)= \alpha_0$$ 
Moreover, by maximum principle and regularity theorem, $v \in C^{\infty}(M)$ 
 and $v>0$. The construction of $\mathcal{F}$ is such that $f$ is critical  
for $\tilde{g}= v^{\frac{4}{N-2}} g$ if and only if $\mathcal{F}(v)$  
is critical for $g$. This is the case here because $\alpha_0$ is critical for  
$g$. Then, the theorem is proved. Thus, in the following, we may assume that 
$v \equiv 0$. This proves step \ref{s2}.\\

\noindent We now assume that $v \equiv 0$ and we prove that 

\begin{step} \label{s3} 
We have :
$$\lim_i \mu_i = K(n,2)^{-2} {(\max_M f )}^{-\frac{2}{N}} \hbox{ and }
\lim_i \int_M v_i^{q_i} dv_g =  
K(n,2)^{-n} {(\max_M f) }^{-\frac{n}{2}}$$
\end{step}

\noindent As easily seen, $\liminf_i \mu_i >0$.
We have, using H\"older inequality:
$$\mu_i=J_i(v_i)=\frac{\int_M {|\nabla v_i|}^2_g dv_g+  t_i \int_M   v_i^2 
dv_g} 
{{\left( \int_M f {v_i}^{q_i} dv_g \right)}^{\frac{2}{q_i}}}$$
$$\geq \frac{\int_M {|\nabla v_i|}^2_g dv_g+  \alpha_0 \int_M   v_i^2 
dv_g} 
{{(\max_M f)}^{\frac{2}{q_i}}
{\left( \int_M  {v_i}^{N} dv_g
  \right)}^{\frac{2}{N}}{Vol(M)}^{1-\frac{q_i}{N} }}+
(t_i-\alpha_0) \frac{     \int_M   v_i^2 
dv_g}   {{\left( \int_M f {v_i}^{q_i} dv_g \right)}^{\frac{2}{q_i}}}   $$ 
Since $\lim_i t_i= \alpha_0$, since $v_i \to 0$ in $L^2(M)$ and since 
$$\liminf_i \int_M fv_i^{q_i} dv_g = \liminf_i
\mu_i^{\frac{q_i}{q_i-2}}>0$$ 
we have
$$\lim_i (t_i-\alpha_0) \frac{     \int_M   v_i^2 
dv_g}   {{\left( \int_M f {v_i}^{q_i} dv_g \right)}^{\frac{2}{q_i}}} =0  $$ 
Moreover,
$$\frac{\int_M {|\nabla v_i|}^2_g dv_g+  \alpha_0 \int_M   v_i^2 
dv_g}
{
{\left( \int_M  {v_i}^{q_i} dv_g
  \right)}^{\frac{2}{q_i}}}=I_{\alpha_0}(v_i) \geq K(n,2)^{-2}$$
because $\alpha_0$ is weakly critical. We obtain that
$$\liminf_i \mu_i \geq  K(n,2)^{-2} {(\max_M f )}^{-\frac{2}{N}} $$
Now, by (\ref{r0}), we can find $w \in
C^{\infty}(M)$ such that 
$$J_{N,\alpha_0}(w) \leq K(n,2)^{-2}  {(\max_M
  f )}^{-\frac{2}{N}} + \epsilon$$ 
where $\epsilon>0$ is as small as wanted.
We have 
$$\limsup_i J_i(w)=J_{N,\alpha_0}(w)  \leq K(n,2)^{-2}  {(\max_M
  f )}^{-\frac{2}{N}} + \epsilon$$
This proves that 
\begin{eqnarray} \label{a0}
\lim_i  \mu_i =  K(n,2)^{-2} {(\max_M f )}^{-\frac{2}{N}}
\end{eqnarray}

\noindent 
Now, we multiply $E(q_i,t_i)$ by $v_i$ and we integrate over $M$. We get:
\begin{eqnarray} \label{a2}
\int_M {|\nabla v_i|}^2_g dv_g + t_i \int_M v_i^2 = \int_M f v_i^{q_i}
dv_g
\end{eqnarray} 
We recall that $\int_M  f v_i^{q_i}
dv_g= \mu_i^{\frac{q_i}{q_i-2}}$. Hence, with H\"older inequality:
 
$$\int_M f v_i^{q_i}
dv_g= \mu_i {\left(\int_M f v_i^{q_i}
dv_g \right)}^{\frac{2}{q_i}} $$
\begin{eqnarray} \label{a3}
\leq  
\mu_i {(\max_Mf)}^{\frac{2}{q_i}} {\left(\int_M v_i^{q_i}
dv_g \right)}^{\frac{2}{q_i}} \leq 
\mu_i {(\max_Mf)}^{\frac{2}{q_i}} {\left(\int_M v_i^{N}
dv_g \right)}^{\frac{2}{N}} Vol(M)^{1-\frac{q_i}{N}}
\end{eqnarray}
 
\noindent Using inequality $S(K(n,2)^2,B_0(g))$ (see introduction), we
obtain that 
$$\int_M f v_i^{q_i}
dv_g $$
$$\leq \mu_i {(\max_Mf)}^{\frac{2}{q_i}} Vol(M)^{1-\frac{q_i}{N}}
\left(K(n,2)^2 \int_M {|\nabla v_i|}^2_g dv_g + B_0(g) \int_M v_i^2 dv_g
\right) $$
Together with (\ref{a2}) and (\ref{a3}), we get 

$$\int_M {|\nabla v_i|}^2_g dv_g + t_i \int_M v_i^2 
\leq \mu_i {(\max_Mf)}^{\frac{2}{q_i}} {\left(\int_M v_i^{q_i}
dv_g \right)}^{\frac{2}{q_i}} $$
\begin{eqnarray} \label{a4}
\leq 
\mu_i {(\max_Mf)}^{\frac{2}{q_i}} Vol(M)^{1-\frac{q_i}{N}}
\left(K(n,2)^2 \int_M {|\nabla v_i|}^2_g dv_g + B_0(g) \int_M v_i^2 dv_g
\right)
\end{eqnarray}
Now, we have 
$J_i(v_i)=\mu_i$ and $\lim_i {\parallel v_i \parallel}_2=0$.  Since
$\int_M f v_i^{q_i} dv_g = \mu_i^{\frac{q_i}{q_i-2}}$, we get from
(\ref{a0}) that  
$$\lim_i \int_M {|\nabla v_i|}^2_g dv_g = \lim_i \mu_i^{\frac{n}{2}}=
{\left( K(n,2)^{-2} {(\max_M f)}^{-\frac{2}{N}} \right)}^{\frac{n}{2}}$$
Taking the limit in $i$ in both sides of inequality (\ref{a4}) and using
(\ref{a0}), 
$${\left( K(n,2)^{-2} {(\max_M f)}^{-\frac{2}{N}} \right)}^{\frac{n}{2}}
\leq K(n,2)^{-2}
 { \left( \lim_i \int_M v_i^{q_i}
  dv_g\right)}^{\frac{2}{N}} $$
$$\leq {\left( K(n,2)^{-2} {(\max_M f)}^{-\frac{2}{N}} \right)}^{\frac{n}{2}} $$
The step then follows immediatly.\\

\noindent Let now $x \in M$ be 
given.  Following usual terminology,  
we say that $x$ is a concentration point if for all $r>0$, 
 
$$\limsup_i \int_{B_x(r)} v_i^{q_i} dv_g > 0$$ 
where $ B_x(r)$ stands for the geodesic ball of center $x$ and radius $r$. 
    
\begin{step} \label{s5} 
Up to a subsequence, $(v_i)$ possesses exactly one concentration point 
$x_0$. Moreover, $x_o$ is a point where $f$ is maximum. If $\bar{\omega} 
\subset
\subset M-\{x_0 \}$ where $\omega$ is an open subset of $M$, 
then $(v_i)$ tends uniformly to $0$ 
with $i$ on $\bar{\omega}$.   
\end{step} 
 
\noindent By step \ref{s3} and 
since $M$ is compact, it is easy  to prove 
the existence of at least one point of 
concentration. We now let $x \in M$ and $r>0$, 
a small positive 
number. Let also $\eta \in C^{\infty}(M)$  a cut-off 
function  supported in $B_x(r)$, such that $0 \leq \eta \leq 1$ and  
$\eta \equiv 1$ on $B_x(\frac{r}{2})$.  
We recall that $v_i$ satisfies equation $E(q_i,t_i)$. We multiply  
$E(q_i,t_i)$ by $\eta^2 v_i^k$ for $k>1$ and integrate over $M$.  
We get: 
\begin{eqnarray} \label{r4} 
\int_M \eta^2 v_i^k \Delta_g v_i dv_g + t \int_M v_i^{k+1} \eta^2 dv_g =  
 \int_M f \eta^2 v_i^{k+q_i-1} dv_g 
\end{eqnarray} 
Integrating by parts, we get  
$$\int_M {|\nabla \eta v_i^{\frac{k+1}{2}} |}_g^2 dv_g = 
\frac{(k+1)^2}{4k} \int_M \eta^2 v_i^k \Delta_g v_i dv_g $$ 
$$+\frac{k+1}{2k} 
\int_M \left( {|\nabla \eta |}_g^2 + \frac{k-1}{k+1} \eta \Delta_g \eta 
\right) v_i^{k+1} dv_g $$ 
    
\noindent Together with (\ref{r4}), this gives 
 
$$\int_M {|\nabla \eta v_i^{\frac{k+1}{2}} |}_g^2 dv_g \leq  
\frac{(k+1)^2}{4k} \int_M f \eta^2 v_i^{k+q_i-1} dv_g$$ 
$$+\frac{k+1}{2k} 
\int_M \left( {|\nabla \eta |}_g^2 + \frac{k-1}{k+1} \eta \Delta_g \eta 
\right) v_i^{k+1} dv_g $$ 
 
\noindent By H\"older inequality, we have  
$$\int_M {|\nabla \eta v_i^{\frac{k+1}{2}} |}_g^2 dv_g \leq  
\frac{(k+1)^2}{4k} \max_M( f)  {\left(\int_M  {(\eta v_i^{\frac{k+1}{2}}
      )}^{q_i}  
dv_g \right)}^{\frac{2}{q_i}}   {\left(\int_{B_x(r)}  v_i^{q_i}  
dv_g \right)}^{\frac{q_i-2}{q_i}}           $$ 
\begin{eqnarray} \label{b1} 
+\frac{k+1}{2k} 
\int_M \left( {|\nabla \eta |}_g^2 + \frac{k-1}{k+1} \eta \Delta_g \eta 
\right) v_i^{k+1} dv_g 
\end{eqnarray}
From inequality $S(K(n,2)^2,B_0(g))$ and again H\"older inequality, we get: 
$$\int_M {|\nabla \eta v_i^{\frac{k+1}{2}} |}_g^2 dv_g \geq 
K(n,2)^{-2}  {\left(\int_M  {(\eta v_i^{\frac{k+1}{2}}
      )}^N  
dv_g \right)}^{\frac{2}{N}} -\alpha_0 \int_M v_i^2 dv_g  $$
 $$\geq K(n,2)^{-2}  {Vol (M)}^{\frac{q_i}{N}-1} {\left(\int_M  {(\eta v_i^{\frac{k+1}{2}}
      )}^{q_i}  
dv_g \right)}^{\frac{2}{q_i}} -\alpha_0 \int_M \eta^2 v_i^{k+1} dv_g $$
Together with (\ref{b1}), we are lead to  

$${ \left(\int_M  {(\eta v_i^{\frac{k+1}{2}}
      )}^{q_i} \right)}^{\frac{2}{q_i}} \left(  K(n,2)^{-2}  {Vol (M)}^{\frac{q_i}{N}-1}-
\frac{(k+1)^2}{4k} \max_{B_x(r)}( f)   {\left(\int_{B_x(r)}  v_i^{q_i}  
dv_g \right)}^{\frac{q_i-2}{q_i}} \right) $$
\begin{eqnarray} \label{r5} 
\leq C \int_M v_i^{k+1} dv_g
\end{eqnarray} 
where $C>0$ is a constant which does not depend on $i$. 
If $x$ is a concentration point, then 
$$\liminf_i \int_{B_x(r)}  v_i^{q_i}  
dv_g >0$$
 Moreover, by step \ref{s3},  
\begin{eqnarray} \label{e1}
\liminf_i \int_{B_x(r)}  v_i^{q_i}  
dv_g \leq  K(n,2)^{-n} {(\max_M f) }^{-\frac{n}{2}}
\end{eqnarray}
 
\noindent Assume that this inequality is strict. Then, if $k$ is
sufficiently close to $1$, we have 
$$\liminf_i \left(  K(n,2)^{-2}  {Vol (M)}^{\frac{q_i}{N}-1}-
\frac{(k+1)^2}{4k} \max_M( f)   {\left(\int_{B_x(r)}  v_i^{q_i}  
dv_g \right)}^{\frac{q_i-2}{q_i}} \right) >0$$
Coming back to (\ref{r5}), we get the existence of $C>0$ independent of $i$
such that 
\begin{eqnarray} \label{r6}
{ \left(\int_M  {(\eta v_i^{\frac{k+1}{2}}
      )}^{q_i} \right)}^{\frac{2}{q_i}} \leq C \int_M v_i^{k+1} dv_g
\end{eqnarray}

\noindent The right hand side of (\ref{r6}) goes to $0$ with $i$. By
H\"older inequality, we would get that 
$$\lim_i \int_{B_x(\frac{r}{2}) }  v_i^{q_i} dv_g \leq 
\lim_i \int_{B_x(\frac{r}{2}) } v_i^{\frac{k+1}{2}q_i} dv_g =0$$
This is impossible since $x$ is a concentration point. It follows that
(\ref{e1}) is a equality and hence, 
there exists one and only one concentration point $x_0$.
Moreover, if $\max_{B_x(r)}( f) < \max_M (f)$ (with $x=x_0$), 
we get (\ref{r6}) in the
same way. Hence, the concentration point $x_0$ is such that 
$\max_M (f)=f(x_0)$.

\noindent Now, let
$\bar{\omega} \subset \subset M-\{x_0\}$ where $\omega$ is an open set of
$M$. Let $0< r < 
dist_g(\omega,x_0)$ and a finite set $(x_j)$ of points of $\omega$ such that 
$$\omega \subset \cup_j B_{x_j}(r)$$
Doing the same with $x=x_j$, this leads to the existence of $C>0$
such that 
$$\int_{\omega} v_i^{\frac{k+1}{2}q_i} dv_g \leq C \int_M v_i^{k+1} dv_g$$
Since $ \frac{k+1}{2}q_i> N+\epsilon$ where $\epsilon >0$ is small, a
simple application of Moser's iterative scheme proves the step.\\

\noindent We now let $s>1$  be a large real number. We claim that 

\begin{step} \label{s6}
For $i$ large enough, the function $\mathcal{F}(v_i^s)$ is weakly critical
for $g$. Moreover, for all $t \in [1,s]$,  $\mathcal{F}(v_i^t) >
\frac{n-1}{4(n-2)}  S_g$.
\end{step}
It is sufficient to prove that for $i$ large enough, $\mathcal{F}(v_i^s) \geq
\alpha_0$. An easy computation gives 
$$\Delta_g (v_i^s)= s v_i^{s-1} \Delta_g v_i - s(s-1) v_i^{s-2} {|\nabla
  v_i|}^2_g \leq s v_i^{s-1} \Delta_g v_i $$
Since $v_i$ satisfies $E(q_i,t_i)$, it follows that 
 
\begin{eqnarray} \label{c1}
\mathcal{F}(v_i^s) \geq s t_i + f(v_i^{\frac{4}{n-2}s}-s v_i^{q_i-2})
\end{eqnarray}

\noindent By step (\ref{s5}), we
know that $v_i$ converges uniformly to $0$ on $\{f \leq 0\}$. Since
$\frac{4}{n-2}s > q_i-2$, we get that, on $\{ f \leq 0\}$ and for $i$ large
enough

\begin{eqnarray} \label{d2} 
\mathcal{F}(v_i^s) \geq s t_i > \alpha_0 
\end{eqnarray}

\noindent Now, we set for $x \geq 0$,
$$\beta(x) = x^{\frac{4}{n-2}s}-s
x^{q_i-2}=x^{q_i-2}(x^{\frac{4}{n-2}s-q_i+2}-s)$$
The minimum of $\beta$ on $[0,+\infty[$ is attained for 
$$x_i= {\left(\frac{(n-2)(q_i-2)}{4}
  \right)}^{\frac{1}{\frac{4}{n-2}s-q_i+2}}$$
Moreover, $x_i \leq 1$ because $q_i \leq N$. Hence,
$|x^i|^{q_i-2} \leq 1$. Hence, 
$$\beta(x_i) \geq  x_i^{\frac{4}{n-2}s-q_i+2}-s= \frac{(n-2)(q_i-2)}{4}-s\geq -s$$
We get from (\ref{c1}) that 
$$\mathcal{F}(v_i^s) \geq s(t_i-\max_M(f))$$
on $\{ f \geq 0\}$.  By (\ref{d1}), $\alpha_0 > \max_M(f)$. Since 
$\lim_i t_i= \alpha_0$, we obtain that  
$$\mathcal{F}(v_i^s) \geq \alpha_0$$
on  $\{ f \geq 0\}$ if $s$ is chosen large enough. By  (\ref{d2}), this
inequality is true on all $M$.
This proves that $\mathcal{F}(v_i^s)$ is weakly critical.

\noindent Now let  
$t \in [1,s]$. 
In the same way, we obtain that if $i$ is large enough, 
$$\mathcal{F}(v_i^t) \geq t (t_i- \max_M(f)) \geq (t_i- \max_M(f))$$
By (\ref{d1}), $\mathcal{F}(v_i^t) > \frac{n-1}{4(n-2)} S_g$.
This proves the step.\\

\begin{step}
Conclusion
\end{step}

\noindent We let $i$ and $s>1$ such that 
$\mathcal{F}(v_i)$ is subcritical, $\mathcal{F}(v_i^s)$ is weakly
critical and for all $t \in [1,s]$, $\mathcal{F}(v_i^t)> \frac{n-1}{4(n-2)}
S_g$. We set $v=v_i$. We also define
$$s_0= \inf \{ t>1 |  \mathcal{F}(v^t) \hbox{ is weakly critical }\}$$
It is clear that
$\mathcal{F}(v^{s_0})$ is weakly critical. Now, let $(t_m)$ a sequence of
real numbers such that 
$1 < t_m < s_0$ and $\lim_{m} t_m=s_0$. We apply theorem 
2.1 with $h_m= \mathcal{F}(v^{t_m})$ and $h=\mathcal{F}(v^{s_0})$.
It follows that $\mathcal{F}(v^{s_0})$ is critical for $g$ and hence that
$f$ is critical for $\tilde{g}=v^{\frac{4}{n-2}s_0} g$. This ends the proof
of main theorem.

%%%%%%%%%%%%%%%%%%%%%%%%%%%%%%%%%%%%%%%%%%%%%%%%%%%%%%%%%%%%%%%%%%%%%%%%%%%%% 

\end{document}